\documentclass[12pt]{amsart}
\usepackage{amssymb,mathrsfs}
\usepackage{amsthm,amsmath,amssymb}
\usepackage[numbers,sort&compress]{natbib}

\makeatletter

\newcommand{\Rmnum}[1]{\expandafter\@slowromancap\romannumeral #1@}
\makeatother
\usepackage{mathtools}

\theoremstyle{plain}
\newtheorem{theorem}{Theorem}[section]
\newtheorem{lemma}[theorem]{Lemma}
\newtheorem{proposition}[theorem]{Proposition}
\newtheorem{corollary}[theorem]{Corollary}

\theoremstyle{definition}
\newtheorem{definition}[theorem]{Definition}

\usepackage[dvips]{graphics}
\usepackage{mathrsfs}
\usepackage{amssymb}
\usepackage{stmaryrd}
\usepackage{pifont}
\usepackage{amsmath}
\usepackage{hyperref}
\usepackage{amsthm}
\usepackage{cases}
\usepackage{mathtools}
\allowdisplaybreaks 
\title[weak L- and M-]{Lattice copies and applications for weak L- and M-weakly compact operators on Banach lattices}

\date{\today}


\author[Z. Wang]{Zhangjun Wang$^{1}$}
\address{$^1$ The first author:School of Mathematics, Southwest Jiaotong University,	Chengdu, Sichuan,	China, 610000.}
\email{zhangjunwang@my.swjtu.edu.cn}

\author[Z. Chen]{Zili Chen$^{2}$}
\address{$^2$ The second author:School of Mathematics, Southwest Jiaotong University, Chengdu, Sichuan,
	China, 611756.}
\email{zlchen@home.swjtu.edu.cn}

\begin{document}

\begin{abstract}
Several recent papers investigated lattice copies and unbounded convergences in Banach lattices. In this paper, we first solve the problem of Rincón‑Villamizar and Leal‑Archila which is an extension of the well-known James distortion theorem. Using lattice copies of $\ell_1$, $\ell_\infty$ and $ba(2^N)$ and unbounded convergence, then we introduce weak L- and M-weakly compact operators on Banach lattices and research the relationship between these operators and L- and M-weakly compact operators. Finally, we study the compactness of weak L-weakly compact and weak M-weakly compact operators.
\end{abstract}
	
\maketitle

\section{Introduction}
Let us begin with some preliminary knowledge to drawing off our research background.

For a set $\Gamma$, $\ell_p(\Gamma) (1\leq p<\infty)$ stands the Banach space of the all bounded families $(a_\gamma)_{\gamma\in\Gamma}$ endowed with the $\ell_p$ norm. When $\Gamma$ is countable, these spaces are denoted as $\ell_p$. If $X$ and $Y$ are Banach spaces, we recall that $Y$ contains a copy 
of $X$ whenever $Y$ contains a subspace isomorphic to $X$; we also say that $Y$ contains almost isometric copies of $X$ if for any $\epsilon>0$, there is an isomorphism $T_\epsilon$ from $X$ into a subspace of $Y$ satisfying $\Vert T_\epsilon\Vert\Vert T^{-1}_\epsilon\Vert\leq1+\epsilon$. Moreover, if $T_\epsilon$ is also lattice isomorphism, $Y$ is said to contains lattice-almost isometric copies of $X$.

For two Banach lattices $E$ and $F$ such that $E$ contains copy of $F$, a problem is determining if $E$ contains lattice-almost isometric copies of $F$. In \cite{J:64}, James proved that if a Banach space $E$ contains a copy of $c_0 (resp. \ell_1)$, then $E$ contains a almost isometric copy of $c_0$ (resp. $\ell_1$). The corresponding statement for $\ell_\infty$ was proved by Partington in \cite{P:81}. In \cite{C:05,C:06}, Chen showed that if a Banach lattice contains a lattice copy of $c_0$ (resp. $\ell_1$, $\ell_\infty$), then it contains lattice-almost isometric copies of $c_0$ (resp. $\ell_1$, $\ell_\infty$). Recently, Rincón‑Villamizar and Leal‑Archila solved the problem (in \cite{X:21}) for $c_0(\Gamma)$, $\ell_1(\Gamma)$ and $\ell_\infty(\Gamma)$ for dual Banach lattice and raise the problem for $\ell_p(\Gamma)$. The aim of Section1 of this paper is the lattice‑almost isometric copies of $\ell_p(\Gamma)$ and solve the problems of Rincón‑Villamizar and Leal‑Archila in \cite{X:21}.


A net $(x_\alpha)$ in a Banach lattice $E$ is \emph{unbounded order (resp. norm, absolute weak)} convergent to some $x$, denoted by $x_\alpha\xrightarrow{uo}x$ (resp. $x_\alpha\xrightarrow{un}x$, $x_\alpha\xrightarrow{uaw}x$), if the net $(|x_\alpha-x|\wedge u)$ converges to zero in order (resp. norm, weak) for all $u\in E_+$. A net $(x_\alpha^\prime)$ in a dual Banach lattice $E^\prime$ is unbounded absolute weak* convergent to some $x^\prime$, denoted by $x_\alpha^\prime\xrightarrow{uaw^*}x^\prime$, if $|x^\prime_\alpha-x^\prime|\wedge u^\prime\xrightarrow{w^*}0$ for all $u^\prime\in E_+^\prime$. For the basic theory of $uo$, $un$, $uaw$ and $uaw^*$-convergence, we refer to \cite{GTX:16,KMT:16,Z:16,T:18}.

In \cite{W:21,W:19}, we studied the continuity functionals and operators for different types of unbounded convergences in Banach lattices, and showed the characterizations of continuous functionals, L-weakly compact sets, L- and M-weakly compact operators on Banach lattices by $uo$, $uaw$ and $uaw^*$-convergence. Based on the above results, we can find that the $un$-convergence is special. Therefore, we will try to study these sets and operators by $un$-convergence. In Section3, we use unbounded convergence and  and lattice copies of $\ell_1$, $\ell_\infty$ and $ba(2^N)$ to introduce and research new classes of sets and operators so called weak (weak*) L-weakly compact sets, weak (weak*) L-weakly compact operators and weak M-weakly compact operators, which contrast with L-weakly compact sets and L-(M-)weakly compact operators. At the end of the paper, we investigate compactness of weak L-weakly compact and weak M-weakly compact operators by these lattice copies. For undefined terminology, notation and basic theory of Riesz space, Banach lattice and linear operator, we refer to \cite{AB:06,MN:91}.

\section{lattice‑almost isometric copies of $\ell_p(\Gamma)$}\label{Sec:2}
Let us determine the lattice‑almost isometric copies of $\ell_p(\Gamma)$ in Banach lattices. If $\tau$ is an ordinal, $|\Gamma|$ (resp. $[\Gamma]^{<\tau}$ ) denotes the cardinality of $\Gamma$ (resp.  the family of all subsets of $\Gamma$ with cardinality less than $\tau$). We denote by $2^\Gamma$ the set of all subsets of $\Gamma$. 
\begin{theorem}\label{1}
	A Banach lattice $E$ contains a lattice copy of $\ell_p(\Gamma)$ ($1\leq p<\infty$) iff it contains lattice-almost isometric copies of $\ell_p(\Gamma)$.
\end{theorem}
\begin{proof}
	Suppose that $E$ contains a lattice copy of $\ell_p(\Gamma)$, then there exists a disjoint family $(x_\gamma)_{\gamma\in\Gamma}$ in $E_+$ and two positive constants $C_1,C_2$ such that $$C_1\sum_{\gamma\in A}|a_\gamma|^p\leq\Vert\sum_{\gamma\in A}a_\gamma x_\gamma\Vert^p\leq C_2\sum_{\gamma\in A}|a_\gamma|^p,$$ for all $A\in2^\Gamma$ and every family of scalars $\{a_\gamma:\gamma\in\Gamma\}$. Let $\tau$ be an ordinal such that $|\Gamma|=\tau$. For $L\in[\Gamma]^{<\tau}$, put $$K_L=\inf\{\Vert\sum_{\gamma\in M}a_\gamma x_\gamma\Vert^p:\sum_{\gamma\in M}a_\gamma^p=1,a_\gamma\geq0, M\in2^\Gamma,L\cap M=\emptyset\}.$$ Clearly, $\{K_L:L\in[\Gamma]^{<\tau}\}$ is bounded. Also, if $L,L^\prime\in [\Gamma]^{<\tau}$ and $L\subset L^\prime$, then $K_L\leq K_{L^\prime}$. Since $C_1\leq K_L\leq C_2$ for all $L$, so we let $$K=\sup\{K_L:L\in[\Gamma]^{<\tau}\}.$$ Take $0<\epsilon<1$ be fixed, and $0<\theta<1<\lambda$ with $\theta/\lambda\geq(1-\epsilon)^p$. Set $L_0\in [\Gamma]^{<\tau}$ such that $\theta K<K_{L_0}\leq K<\lambda K$. Thus there is $M_0\in2^\Gamma$ with $M_0 \cap L_0 =\emptyset$ and $(b_\gamma^{M_0})_{\gamma\in M_0}$ such that $$\theta K<K_{L_0}\leq\Vert\sum_{\gamma\in M_0}b_\gamma^{M_0} x_\gamma\Vert^p<\lambda K, \sum_{\gamma\in M_0}(b_\gamma^{M_0})^p=1.$$ Now let $\alpha-1=\beta$ be an ordinal with $\alpha<\tau$. Assmue that $\{M_\eta\in 2^\Gamma:\eta\leq\beta\}$ has been defined in a way such that:
	\begin{enumerate}
		\item for each $\eta\leq\beta$, we have $$\theta K<K_{L_0}\leq\Vert\sum_{\gamma\in M_\eta}b_\gamma ^{M_\eta}x_\gamma\Vert^p<\lambda K, \sum_{\gamma\in M_\eta}(b_\gamma^{M_\eta})^p=1.$$
		\item The family $\{M_\eta\in 2^\Gamma:\eta\leq\beta\}$ is a family of disjoint finite sets and $M_\eta\cap L_0=\emptyset$ for each $\eta\leq\beta$.
	\end{enumerate}
	Since $|\cup_{\eta\leq\beta}M_\eta|<|\Gamma|$, so $\Gamma/\cup_{\eta\leq\beta}M_\eta\neq\emptyset$. If $N_\beta=\cup_{\eta\leq\beta}M_\eta$, then $\theta K<K_{L_0}\leq K_{N_\beta\cup L_0}<\lambda K$. Hence, there exists $M\in 2^\Gamma$ such that $M\cap( N_\beta\cup L_0)=\emptyset$ and $$\theta K< K_{L_0}\leq K_{N_\beta\cup L_0}\leq \Vert \sum_{\gamma\in M}b_\gamma^Mx_\gamma\Vert^p<\lambda K, \sum_{\gamma\in M}(b_\gamma^{M})^p=1.$$ Therefore, we take $M_\alpha:=M$. Now, we suppose that $\alpha$ is a limit ordinal and $\{M_\beta:\beta<\alpha\}$ has been defined. It follows from $$|\cup_{\beta<\alpha}M_\beta|\leq\alpha<\tau=|\Gamma|$$ that $\Gamma/\cup_{\beta<\alpha}M_\beta\neq\emptyset$. If $N_\alpha=\cup_{\beta<\alpha}M_\beta$, then $\theta K<K_{L_0}\leq K_{N_\alpha\cup L_0}<\lambda K$. So, there is $M^\prime\in2^\Gamma$ such that $M^\prime\cap(N_\alpha\cup L_0)=\emptyset$ and $$\theta K< K_{L_0}\leq K_{N_\beta\cup L_0}\leq \Vert \sum_{\gamma\in M^\prime}b_\gamma^{M^\prime}x_\gamma\Vert^p<\lambda K, \sum_{\gamma\in M^\prime}(b_\gamma^{M^\prime})^p=1.$$Let $M_\alpha:=M^\prime$. By the way, we have constructed a family $\mathbb{M}=\{M_\alpha:\alpha<\tau\}=\{M_\gamma:\gamma\in\Gamma\}$ in $2^\Gamma$ and a family $\mathbb{Y}=\{y_{M_\gamma}:\gamma\in\Gamma\}$ in $E_+$ with $|\mathbb{M}|=|\mathbb{Y}|=|\Gamma|$ satisfying 
	\begin{enumerate}
		\item for any $\gamma\in\Gamma$, there exists a set of positive scalars $\{b^{M_\gamma}_{\gamma^\prime}:\gamma^\prime\in M_\gamma\}$ for all $\gamma^\prime\in M_\gamma$ such that
		$$y_{M_\gamma}=\sum_{\gamma^\prime\in M_\gamma}b^{M_\gamma}_{\gamma^\prime}x_{\gamma^\prime}, \sum_{\gamma^\prime\in M_\gamma}(b^{M_\gamma}_{\gamma^\prime})^p=1.$$ and $\theta K<K_{L_0}\leq\Vert y_{M_\gamma}\Vert^p<\lambda K$.
		\item The family $\mathbb{M}$ is a family of disjoint finite sets with $M_\gamma\cap L_0=\emptyset$ for each $\gamma\in\Gamma$.
	\end{enumerate}
	Let $z_\gamma=(\lambda K)^{-1/p}y_{M_\gamma}$, clearly,  $(z_\gamma:\gamma\in\Gamma)$ is a disjoint family in $B_E^+$. For each $\mu:=(\mu_\gamma)_{\gamma\in\Gamma}\in l_p(\Gamma)$, we have $$\Vert\sum_{\gamma\in\Gamma}\mu_\gamma z_\gamma\Vert^p\leq \Vert\sum_{\gamma\in\Gamma}\mu_\gamma\Vert^p(\sup_{\gamma\in\Gamma}\Vert z_\gamma\Vert^p)\leq \Vert\sum_{\gamma\in\Gamma}\mu_\gamma\Vert^p.$$
	On the other hand,
	\begin{align*}
		\Vert\sum_{\gamma\in\Gamma}\mu_\gamma z_\gamma\Vert^p&=(\lambda K)^{-1}	\Vert\sum_{\gamma\in\Gamma}\mu_\gamma y_\gamma\Vert^p\\
		&=(\lambda K)^{-1}\Vert\mu\Vert^p\cdot\big\Vert\sum_{\gamma\in\Gamma}|\mu_\gamma|y_\gamma/\Vert\mu\Vert\big\Vert^p\\
		&\geq(\lambda K)^{-1}K_{L_0}\Vert\mu\Vert^p\geq(1-\epsilon)^p\Vert\mu\Vert^p.
	\end{align*}
	Now we define $T_\epsilon\mu=\sum_{\gamma\in\Gamma}\mu_\gamma z_\gamma$, it can be easily verified that, $T_\epsilon$ is a lattice embedding from $\ell_p(\Gamma)$ into $E$ such that $(1-\epsilon)\Vert\mu\Vert\leq\Vert T_\epsilon\mu\Vert\leq\Vert\mu\Vert$. The proof is completed.
\end{proof}	

\section{weak L-weakly compact and weak M-weakly compact operators}\label{Subsec:1}

A vector $e$ in a Riesz space $E$ is said to be \emph{strong order unit} if the ideal $I_e$ generated by $e$ equal to $E$. According to \cite[Theorem~4.21 and Theorem~4.29]{AB:06}, a Banach lattice $E$ is lattice and norm isomorphic to $C(K)$ for some compact Hausdorff space $K$ whenever $E$ has a strong order unit. It follows from \cite[Theorem~2.3]{KMT:16} that the $un$-convergence implies norm convergence in $E$ iff $E$ has strong order unit. 

Since every bounded $un$-null sequence in $l_\infty$ is norm convergent to zero. It is clearly that every $un$-null sequence in $B_{l_\infty}$ converges uniformly to zero on $B_{l_1}$ and $B_{ba(2^N)}$. So for the identical operator $I:l_1\rightarrow l_1$, every $un$-null sequence in $B_{l_\infty}$ converges uniformly to zero on $I(B_{l_1})$ and $I^{\prime\prime}(B_{ba(2^N)})$. And $I^\prime x_n^\prime\rightarrow0$ for every $un$-null sequence $(x_n^\prime)\subset B_{l_\infty}$. 

Clearly, $B_{l_1}$ and $B_{ba(2^N)}$ are not L-weakly compact sets, $I$ and $I^{\prime\prime}$ are not L-weakly compact operators and $I^\prime$ is not M-weakly compact. Therefore, we introduce these sets and operators.
\begin{definition}
	Let $E$ be a Banach lattice, a bounded subset $A\subset E (B\subset E^\prime)$ is called \emph{weak(weak*) L-weakly compact} set whenever $\sup_{x\in A}|x^\prime_n(x)|\rightarrow0 (\sup_{x^\prime\in B}|x^\prime(x_n)|\rightarrow0)$ for every bounded $un$-null sequence $\{x_n^\prime\}\big(\{x_n\}\big)$ in $E^\prime (E)$. Clearly, $A (B)$ is weak (weak*) L-weakly compact set iff for each sequence $\{x_n\}\big(\{x_n^\prime\}\big)$ in $A (B)$, $x_n^\prime(x_n)\rightarrow0$ for every bounded $un$-null sequence $\{x_n^\prime\}\big(\{x_n\}\big)$ in $E^\prime (E)$.
	
	Respectively, let $X$ be a Banach space. A continuous operator $T:X\rightarrow E (T:X\rightarrow E^\prime)$ is said to be \emph{weak (weak*) L-weakly compact } operator if $T(B_X)$ is weak (weak*) L-weakly compact set in $E (E^\prime)$. Clearly, $T$ is weak (weak*) L-weakly compact iff $f_n(Tx_n)\rightarrow0$ for every sequence $(x_n)\subset B_X$ and every bounded $un$-null sequence $(f_n)$ in $E^\prime (E)$.
	
	A continuous operator $T:E\rightarrow Y$ is said to be \emph{weak M-weakly compact operator} whenever $Tx_n\rightarrow0$ for every $un$-null sequence $(x_n)\subset B_E$.
\end{definition}
For an operator $T : E \rightarrow F$ between two Riesz spaces
we shall say that its modulus $|T|$ exists (or that $T$ possesses a modulus)
whenever $|T| := T \vee (-T)$ exists. The carrier of $T$ is denoted by $C_T$ with $C_T:=\{x\in E:|T|(|x|)=0\}^d$. Accroding to \cite[Theorem~2.2]{W:21}, the carriers of the $uo$-continuous, $un$-continuous, $uaw$-continuous, $uaw^*$-continuous and disjoint continuous operators on atomic Banach lattice are finite-dimensional. Hence, we assume that these unbounded convergence sequences for these operators are norm bounded.

The following basic properties of weak and weak* L-weakly compact sets in Banach lattice can be obtained.

\begin{proposition}\label{solid-convex hull}
	For a Banach lattice $E$, the following statements hold.
	\begin{enumerate}
		\item Every subset and finite union of weak L-weakly compact set in $E$ is weak L-weakly compact set in $E$.
		\item Every subset and finite union of weak* L-weakly compact set in $E^\prime$ is weak* L-weakly compact set in $E^\prime$.
		\item The solid convex hull of weak L-weakly compact set in $E$ is weak L-weakly compact set in $E$.
		\item The solid convex hull of weak* L-weakly compact set in $E^\prime$ is weak* L-weakly compact set in $E^\prime$.
	\end{enumerate}
\end{proposition}
\begin{proof}
	$(1)$ and $(2)$. Obvious.
	
	$(3)$ For a weak L-weakly compact subset $A$ of Banach lattice $E$, let $Sol(A)$ denote the solid hull of $A$.  $\sup_{x\in Sol(A)}|x_n^\prime(x)|=\sup_{x\in A}\{|y_n^\prime(x)|:|y_n^\prime|\leq |x_n^\prime|\}$. Clearly, $y_n^\prime\xrightarrow{un}0$ since $x_n^\prime\xrightarrow{un}0$. So we have $\sup_{x\in Sol(A)}|x_n^\prime(x)|\rightarrow0$, therefore $Sol(A)$ is also weak L-weakly compact set. 
	
	let $B$ denote the solid convex hull of $A$ as
	$$B:=\{\sum_{i=1}^{n}\lambda_i x_i:x_i\in Sol(A),\lambda_i\geq0,\sum_{i=1}^{n}\lambda_i=1\}.$$
	Then $|x^\prime_n(x)|=|x^\prime_n(\sum_{i=1}^{n}\lambda_i x_i)|\leq \sum_{i=1}^{n}|\lambda_i x_n^\prime(x_i)|\rightarrow0$ for every bounded $un$-null sequence $(x_n^\prime)$ in $E^\prime$. Hence $\sup_{x\in B}|x^\prime_n(x)|\rightarrow0$, therefore $B$ is weak L-weakly compact set.
	
	$(4)$ is similar to $(3)$. 
\end{proof}

According to the above results, it is clear that every L-weakly compact set and operator is weak and weak* L-weakly compact and every M-weakly compact operator is weakly M-weakly compact, but the converse does not hold in general. Then, we consider that when are weak (weak*) L- and M-weakly compact L- and M-weakly compact.

The following results are some characterizations of order continuous Banach lattice by weak L-weakly compact, weak* L-weakly compact and weak M-weakly compact operators.

\begin{lemma}\label{weak L-strong continuous dual}
	Let $E$ and $F$ be Banach lattices, the following holds.
	\begin{enumerate}
		\item A continuous operator $T:E\rightarrow F$ is weak L-weakly compact iff $T^\prime:F^\prime\rightarrow E^\prime$ is weak M-weakly compact.
		\item A continuous operator $T:E\rightarrow F$ is weak M-weakly compact iff $T^\prime:F^\prime\rightarrow E^\prime$ is weak* L-weakly compact.
	\end{enumerate}
\end{lemma}
\begin{proof}
	$(1)$ Since $y^\prime(y)=y^\prime(Tx)=(T^\prime y^\prime)(x)$, hence $\Vert T^\prime y_n^\prime\Vert=\sup_{y\in T(B_E)}|y_n^\prime(y)|$. Assume that $T$ is weak L-weakly compact, then for a $un$-null sequence $(y_n^\prime)\subset B_{F^\prime}$, we have $\Vert T^\prime y_n^\prime\Vert=\sup_{y\in T(B_E)}|y_n^\prime(y)|\rightarrow0$, hence $T^\prime$ is weak M-weakly compact. The converse is similar. 
	
	$(2)$ is similar to $(1)$.
\end{proof}
\begin{theorem}\label{strong-M-weakly}
	Let $E$ be a Dedekind $\sigma$-complete Banach lattice, the following conditions are equivalent.
	\begin{enumerate}
		\item $E$ has order continuous norm.
		\item For each Banach space $F$, every weak M-weakly compact operator $T:E\rightarrow F$ is M-weakly compact.
		\item Every positive weak M-weakly compact operator $T:E\rightarrow l_\infty$ is M-weakly compact.
		\item For each Banach space $X$, every adjoint weak* L-weakly compact operator $T^\prime:X^\prime\rightarrow E^\prime$ for continuous operator $T:E\rightarrow X$ is L-weakly compact.
		\item Every positive adjoint weak* L-weakly compact operator $T^\prime:ba(2^N)\rightarrow E^\prime$ for continuous operator $T:E\rightarrow l_\infty$ is L-weakly compact.
	\end{enumerate}
\end{theorem}
\begin{proof}
	$(1)\Rightarrow(2)$ For a weak M-weakly compact operator $T:E\rightarrow F$ and a bounded disjoint sequence $(x_n)\subset E$. It is easy to see that $x_n\xrightarrow{un}0$ by \cite[Proposition~3.5]{KMT:16}. So we have $Tx_n\rightarrow0$, therefore $T$ is M-weakly compact.
	
	$(2)\Rightarrow(3)$ Obvious.
	
	$(3)\Rightarrow(1)$ Assume that $E$ is not order continuous, it follows from \cite[Theorem~4.51]{AB:06} that $E$ contains lattice copy of $l_\infty$. According to \cite[Proposition~1.5.10(1)]{MN:91}, the identical operator $I:l_\infty\rightarrow l_\infty$ can extension to all of $E$. Moreover, $I$ has a positive extension to all of the $E$ by \cite[Exercise~1.5.E1]{MN:91}. Therefore, there exists a positive projection $P:E\rightarrow l_\infty$. 
	
	Let $T=P:E\rightarrow l_\infty$. For a bounded $un$-null sequence $(x_n)\subset E$, it follows from \cite[Theorem~4.3]{KMT:16} that $Tx_n\xrightarrow{un}0$. Since $l_\infty$ has strong order unit, hence $Tx_n\rightarrow0$ by \cite[Theorem~2.3]{KMT:16}. Hence, $T$ is positive weak M-weakly compact. But, $T$ is not M-weakly comopact. Indeed, for the disjoint unit vectors $(e_n)$ of $l_\infty$, $\Vert Te_n\Vert=1\nrightarrow0$, this leads to contradiction. Therefore, $E$ has order continuous norm.
	
	$(1)\Rightarrow(4)$. Every disjoint sequence $(x_n)\subset E$ is $un$-null.
	
	$(4)\Rightarrow(5)$. Obvious.
	
	$(5)\Rightarrow(4)$ According to Lemma~\ref{weak L-strong continuous dual}.
\end{proof}

Dually, we have a similar result of dual Banach lattice.

\begin{theorem}\label{weak-L-weakly compact operator}
	Let $E$ be a Banach lattice, then the following statements are equivalent.
	\begin{enumerate}
		\item $E^\prime$ has order continuous norm.
		\item For each Banach space $X$, every weak L-weakly compact operator $T:X\rightarrow E$ is L-weakly compact.
		\item Every positive weak L-weakly compact operator $T:l_1\rightarrow E$ is L-weakly compact.
	\end{enumerate}
\end{theorem}
\begin{proof}
	$(1)\Rightarrow(2)$. Every disjoint sequence $(x_n^\prime)\subset E^\prime$ is $un$-null.
	
	$(2)\Rightarrow(3)$. Obvious.
	
	$(3)\Rightarrow(1)$. For a positive weak L-weakly compact operator $T$, according to Lemma\ref{weak L-strong continuous dual}, $T^\prime:E^\prime\rightarrow l_\infty$ is positive weakly M-weakly compact. It follows from Theorem\ref{strong-M-weakly} that $T^\prime$ is M-weakly compact. Therefore, $T$ is L-weakly compact by \cite[Theorem~5.64]{AB:06}
\end{proof}

The following results are some characterizations of weak (weak*) L-weakly compact sets and weak (weak*) L-weakly compact operators about disjoint sequence.

\begin{theorem}\label{disjoint-weak L}
	Let $E$ be a Banach lattice, $A$ a bounded solid subset of $E$ and $B$ a bounded solid subset of $E^\prime$. The following statements hold.
	\begin{enumerate}
		\item If $E$ has order continuous norm, then $A$ is weak L-weakly compact iff $x_n^\prime(x_n)\rightarrow0$ for every positive disjoint sequence $(x_n)$ in $A$ and each bounded $un$-null sequence $(x_n^\prime)$ in $E^\prime$.
		
		\item If $E^\prime$ has order continuous norm, then $B$ is weak* L-weakly compact iff $x_n^\prime(x_n)\rightarrow0$ for every positive disjoint sequence $(x_n^\prime)$ in $B$ and each bounded $un$-null sequence $(x_n)$ in $E$
	\end{enumerate}
\end{theorem}

\begin{proof}
	$(1)\Rightarrow$. Clearly.
	
	$(1)\Leftarrow$ Let $(x_n^\prime)$ be a bounded $un$-null sequence in $E^\prime$. To finish the proof, we have to show that  $\sup_{x\in A}|x_n^\prime(x)|\rightarrow0$. Assume by way of contradiction that $\sup_{x\in A}|x_n^\prime(x)|\nrightarrow0$.  Then, by passing to a subsequence if necessary, we can suppose that there would exist some $\epsilon>0$ such that $\sup_{x\in A}|x_n^\prime(x)|>\epsilon$ for all $n$. Note that the equality $\sup_{x\in A}|x_n^\prime(x)|=\sup_{0\leq x\in A}|x_n^\prime|(x)$ holds, since $A$ is solid. $|x_n^\prime|\xrightarrow{w^*}0$ since $E^\prime$ is order continuous. Let $n_1=1$. Because $|x_n^\prime|(4x_{n_1})\rightarrow0$, there exists some $1<n_2\in\mathbb{N}$ such that $|x_{n_2}^\prime|(4x_{n_1})<\frac{1}{2}$. It is easy to see that we can find a strictly increasing subsequence $(n_k)_{k=1}^\infty\subset \mathbb{N}$ such that $|x_{n_{m+1}}^\prime|(4^m\sum_{k=1}^{m}x_{n_k})<\frac{1}{m}$ for all $m$. Let
	$$x=\sum\limits_{k=1}^\infty2^{-k}x_{n_{k}}, y_m=(x_{n_{m+1}}-4^m\sum\limits_{k=1}^{m}x_{n_k}-2^{-m}x)^+.$$
	According to \cite[Lemma~4.35]{AB:06},$(y_m)$ is a disjoint sequence in $A\cap E_+$. Now, we have
	\begin{align*}
		|x^\prime_{n_{m+1}}|(y_m)&=|x^\prime_{n_{m+1}}|(x_{n_{m+1}}-4^m\sum\limits_{k=1}^{m}x_{n_k}-2^{-m}x)^+\\
		&\geq |x^\prime_{n_{m+1}}|(x_{n_{m+1}}-4^m\sum\limits_{k=1}^{m}x_{n_k}-2^{-m}x)\\
		&=|x^\prime_{n_{m+1}}|(x_{n_{m+1}} )-|x^\prime_{n_{m+1}}|(4^m\sum\limits_{k=1}^{m}x_{n_k})-2^{-m}|x^\prime_{n_{m+1}}|x\\
		&>\epsilon-\frac{1}{m}-2^{-m}|x^\prime_{n_{m+1}}|x.
	\end{align*} 
	Let $m\rightarrow\infty$, it is clear that $2^{-m}|x^\prime_{n_{m+1}}|x\rightarrow0$. Hence, $|x^\prime_{n_{m+1}}|(y_m)\nrightarrow0$. This leads to a contradiction.
	
	The proof of $(2)$ is similar.
\end{proof}

Using similar proof methods, we also have the following result.
\begin{theorem}\label{weak L-weakly compact operator-disjoint}
	Let $E$ and $F$ be Banach lattices, for a positive operator $T:E\rightarrow F$, the following statements hold.
	\begin{enumerate}
		\item If $F$ has order continuous norm, then $T$ is weak L-weakly compact iff $y_n^\prime(Tx_n)\rightarrow0$ for each positive disjoint sequence $(x_n)\subset B_E$ and every bounded $un$-null sequence $(y_n^\prime)$ in $F^\prime$.
		\item For the positive adjoint operator $T^\prime:F^\prime\rightarrow E^\prime$ of $T$, if $E^\prime$ has order continuous norm, then $T$ is weak* L-weakly compact iff $T^\prime y_n^\prime(x_n)\rightarrow0$ for each positive disjoint sequence $(y_n^\prime)\subset B_{F^\prime}$ and every bounded $un$-null sequence $(x_n)$ in $E$.
	\end{enumerate}
\end{theorem}
\begin{proof}
	$(1)\Rightarrow$ Obvious.
	
	$(1)\Leftarrow$ Let $(y_n^\prime)$ be an arbitrary bounded $un$-null sequence in $F^\prime$. $|y_n^\prime|\xrightarrow{w^*}0$ since $F$ is order continuous. Hence, $|T^\prime(y_n^\prime)|(z)=\sup_{y\in T[-z,z]}|y_n^\prime(y)|\rightarrow0$ for each $z\in E_+$. Without loss of generality, $y_n^\prime\geq0$ for all $n$. To finish the proof, we have to show that $\sup_{x\in B_E}|y_n^\prime(Tx)|\rightarrow0$. Assume by way of contradiction that $\sup_{x\in B_E}|y_n^\prime(Tx)|\nrightarrow0$. Then, by passing to a subsequence if necessary, we can suppose that there would exist some $\epsilon>0$ such that $\sup_{x\in B_E}|y_n^\prime(Tx)|>\epsilon$ for all $n$. Note that the equality  $\sup_{x\in B_E}|y_n^\prime(Tx)|=\sup_{0\leq x\in B_E}\{|y^\prime_n\circ T|(x)\}$ since $A$ is solid. For every $n$, there exists $z_n$ in $B_E\cap E_+$ such that $|T^\prime( y^\prime_n)|(z_n)>\epsilon$. It is similar to the proof of Theorem \ref{disjoint-weak L} that there exists a subsequence $(y_n)$ of $(z_n)$ and a subsequence $(g_n)$ of $( y_n^\prime)$ such that 
	$$|g_n\circ T |(y_n) > \epsilon, |g_{n+1} \circ T |(4^n\sum_{i=1}^{n}y_i)<\frac{1}{n}.$$
	Let $x=\sum_{i=1}^{\infty}2^{-i}y_i$ and $x_n=(y_{n+1}-4^n(\sum_{i=1}^{n}y_i)-2^{-n}x)^+$, according to \cite[Lemma~4.35]{AB:06}, $(x_n)$ is positive and disjoint. Hence, 	$$|g_{n+1} \circ T |(x_n)\geq|g_{n+1} \circ T |(y_{n+1}-4^n(\sum_{i=1}^{n}y_i)-2^{-n}x)> \epsilon-\frac{1}{n}-2^{-n}|g_{n+1} \circ T |x.$$
	Therefore, $|g_{n+1} \circ T|(x_n)\nrightarrow0$. Clearly, there exists a sequence $(u_n)$ in $E$ satisfying $|u_n|\leq x_n$ such that $|g_{n+1}(Tu_n)|=|g_{n+1} \circ T|(x_n)$. As applications of $$|g_{n+1}(Tu_n^+)|+|g_{n+1}(Tu_n^-)|\geq|g_{n+1}(Tu_n)|=|g_{n+1} \circ T|(x_n)\nrightarrow0,$$ we have $g_{n+1}(Tu_n^+)\nrightarrow0$. This leads to a contradiction.
	
	The rest of the proof is similar.
\end{proof}

\section{compactness of weak L-weakly compact and weak M-weakly compact operators}\label{}

Compact operator is not weak L-weakly compact and weak M-weakly compact in general. Considering the compact operator (rank is 1) $T:l_1\rightarrow l_\infty$ define as $T(x_n)=(\sum_{n=1}^{\infty}x_n,\sum_{n=1}^{\infty}x_n,...)\in l_\infty$ for each $(x_n)\in l_1$. It is clear that $T$ is not weak (weak*) L-weakly compact and weak M-weakly compact.

Weak (weak*) L-weakly compact and weak M-weakly compact operators are also not compact in general. Considering the identical operator $I:l_1\rightarrow l_1$, $I$ is weak L-weakly compact, $I^\prime$ is weak M-weakly compact and $I^{\prime\prime}$ is weak* L-weakly compact. But $I$, $I^\prime$ and $I^{\prime\prime}$ are not compact.

In this section, we research when weak (weak*) L-weakly compact and weak M-weakly compact operator is compact and the converse.

It is easy to see that every semi-compact operator is L-weakly compact whenever the range space is order continuous. Therefore, the following result can be obtained immediately.

\begin{proposition}
	Let $X$ be a Banach space and $E$ be a Banach lattice, then the following hold.
	\begin{enumerate}
		\item If $E$ has order continuous norm, then every compact operator from $X$ into $E$ is weak L-weakly compact.
		\item If $E^\prime$ has order continuous norm, then every compact operator from $X$ into $E^\prime$ is weak* L-weakly compact.
		\item If $E^\prime$ has order continuous norm, then every compact operator from $E$ into $X$ is weak M-weakly compact.
	\end{enumerate}
\end{proposition}

Recall that a vector $e > 0$ in an Banach lattice lattice $E$ is an atom if for any
$u, v \in [0, e] $ with $u \wedge v = 0$, either $u = 0$ or $v = 0$. In this case, the band generated by $e$ is $span\{e\}$. Moreover, the band projection $P_e:E\rightarrow span\{e\}$ defined by
$$P_ex = \sup_n(x^+\wedge ne)-\sup_n(x^-\wedge ne)$$
exists, and there is a unique positive linear functional $f_e$ on $E$ such that $P_e(x) = f_e(x)e$
for all $x\in E$. We call $f_e$ the coordinate functional with the atom $e$. Clearly, the span of any finite set of atoms is also a projection band. A Banach lattice $E$ is called \emph{atomic} if $E$ has a complete disjoint system consisting of atoms. 

The order continuous part $E^a$ of a Banach lattice $E$ is given by $$E^a=\{x\in E:\text{every monotone increasing sequence in $[0,|x|]$ is norm convergent}\}.$$ According to \cite[Corollary~2.3.6]{MN:91}, it is equivalent to $$E^a=\{x\in E:\text{every disjoint sequence in $[0,|x|]$ is norm convergent}\}.$$ A Banach lattice $E$ is said to be \emph{order continuous} whenever $\Vert x_\alpha\Vert\rightarrow0$ for every net $x_\alpha\downarrow0$ in $E$. By \cite[Proposition~2.4.10]{MN:91}, $E^a$ is the largest closed ideal with order continuous norm of $E$.

It is natural to consider that when are weak(weak*) L-weakly and M-weakly compact operator compact. The following results answer the question.

\begin{theorem}\label{wl is compact}
	For a Banach lattice $F$, the following statements are equivalent.
	\begin{enumerate}
		\item $F^a$ is atomic and $F^\prime$ has order continuous norm.
		\item For each Banach space $X$, every weak L-weakly compact operator $T:X\rightarrow F$ is compact.
		\item For each Banach lattice $E$ without order continuous dual, every positive weak L-weakly compact operator $T:E\rightarrow F$ is compact.
	\end{enumerate}
\end{theorem}
\begin{proof}
	$(1)\Rightarrow (2)$ Since $F^\prime$ has order continuous norm, so every weak L-weakly compact operator $T:X\rightarrow F$ is L-weakly compact by Theorem \ref{weak-L-weakly compact operator}.  According to \cite[Proposition~3.6.2]{MN:91}, for any $\epsilon>0$, there exists some $y\in F^a_+$ such that $$T(B_E)\subset [-y,y]+\epsilon\cdot B_E.$$ Since $F^a$ is atomic, it follows from \cite[Theorem~6.1(5)]{W:99} that the order interval $[-y,y]$ is norm compact. Using \cite[Theorem~3.1]{AB:06}, we have $T(B_E)$ is relatively compact set in $F$, so $T$ is also compact.
	
	$(2)\Rightarrow (3)$. Obvious.
	
	$(3)\Rightarrow (1)$ We claim that $F^\prime$ is order continuous. Assume that $F^\prime$ is not order continuous, then $E^\prime$ and $F^\prime$ contain the lattice copy of $l_1$, moreover there exists a positive projection from $E^\prime$ to $l_1$. Let $R:E\rightarrow l_1$ be the positive projection, $S:l_1\rightarrow F$ be the canonical injection from $l_1$ into $F$ and $T=S\circ R:E\rightarrow l_1\rightarrow F$. It is clear that $T$ is weak L-weakly compact since $T^\prime$ is weak M-weakly compact, but not compact. Therefore, $F^\prime$ has order continuous norm.
	
	Then we prove that $F^a$ is atomic. Since $E^\prime$
	is not order continuous then, by \cite[Theorem ~2.4.14]{MN:91}, there is a norm bounded disjoint
	sequence $(u_n)$ of positive elements in $E$ which does not converge weakly to zero. Without loss of generality, we may assume that $\Vert u_n\Vert\leq 1$ for all $n$ and that there are $\phi\in E^\prime_+$ and $\epsilon>0$ such that $\phi (u_n) > \epsilon$ for all $n$. It follows from \cite[Theorem~116.3]{LZ:83} that the components $\phi_n$ of $\phi$, in the carriers $C_{u_n}$, form an order bounded
	disjoint sequence in $E^\prime_+$ such that \text{$\phi_n (u_n) = \phi (u_n)$ for all $n$ and $\phi_n (u_m) = 0$ if $n\neq m$.} Note that $0 \leq \phi_n \leq \phi$ for all $n$. 
	
	Assume that $F^a$ is not atomic, it follows from \cite[Theorem~6.1]{W:99} that there
	exists some $0 \leq y \in F^a$
	such that $[0, y]$ is not norm compact. Now, fix a sequence $(y_n)$
	in $[0, y]$ which has no norm convergent subsequence in $F^a$ and none in $F$.
	
	Define an operator $T : E \rightarrow F$ by
	$$T(x)=\sum_{n=1}^{\infty}(\frac{\phi_n(x)}{\phi(u_n)})y_n$$
	for $x\in E$. Note that in view of the inequality
	$$\sum_{n=1}^{\infty}\Vert (\frac{\phi_n(x)}{\phi(u_n)})y_n\Vert\leq \frac{1}{\epsilon}\Vert y\Vert\sum_{n=1}^{\infty}\phi_n(|x|)\leq \frac{1}{\epsilon}\Vert y\Vert\phi(|x|)$$
	for each $x \in E,$ the series defining $T$ converges in norm for each $x \in E$ and $T (u_n) = y_n$ for all $n$. Hence the operator $T$ is well defined and it is also easy to see that $T$ is a
	positive operator.
	
	Since $(y_n)$ has no norm convergent subsequence in $F$, so $T$ is not compact by Grothendieck theorem (\cite[Theorem~5.3]{AB:06}). However, $T$ maps norm bounded subset in $E$ to an order bounded subset in $F$. To see this, note that for all $x \in B_E$, we have
	$$|T(x)| \leq T(|x|) =\sum_{n=1}^{\infty}(\frac{\phi_n(|x|)}{\phi(u_n)})y_n\leq \frac{1}{\epsilon}\big(\sum_{n=1}^{\infty}\phi_n(|x|)\big)y\leq \frac{1}{\epsilon}\phi(|x|)y\leq \frac{1}{\epsilon}\Vert\phi\Vert y.$$
	So $T$ is L-weakly compact operator, hence $T$ is weak L-weakly compact operator. This leads contradiction, so $F^a$ is atomic.
\end{proof}

The following result shows that when do weakly L-weakly compact operators and compact operators coincide.
\begin{theorem}\label{compact=w L-weakly compact}
	Let $F$ be a Banach lattice, then the following statements are equivalent.
	\begin{enumerate}
		\item $F$ is atomic and both $F$ and $F^\prime$ are order continuous.
		\item For each Banach space $X$, every continuous operator $T:X\rightarrow F$ is compact operator iff $T$ is weak L-weakly compact.
		\item For each Banach lattice $E$ without order continuous dual, every positive operator $T:E\rightarrow F$ is compact operator iff $T$ is weak L-weakly compact.
	\end{enumerate}
\end{theorem}

\begin{proof}
	$(1)\Rightarrow(2)$. Since $F$ is order continuous, so every compact operator $T:E\rightarrow F$ is weak L-weakly compact. Since $F$ is an atomic Banach lattice with order continuous dual, it follows form Theorem \ref{wl is compact} that every weak L-weakly compact operator is compact.
	
	$(2)\Rightarrow(3)$. Obvious.
	
	$(3)\Rightarrow(1)$. According to Theorem \ref{wl is compact}, we have $F^a$ is atomic and $F^\prime$ is order continuous. We claim that $F$ has order continuous norm. 
	
	Suppose that $F$ is not order continuous, by \cite[Corollary~2.4.3]{MN:91}, there exists a disjoint sequence $(x_n^\prime)\subset B_{F^\prime}$ such that $x_n^\prime\nrightarrow0$ in $\sigma(F^\prime,F)$. That is, there is some $y\in F_+$ with $x^\prime_n(y)\nrightarrow0$. As $E\ne\{0\}$, we may fix
	$u\in E$ and pick a $\phi\in E^\prime$
	such that $\phi(u)=\Vert u\Vert=1$ holds.
	
	Now, we consider operator $T : E\rightarrow F$ defined by $$T(x)=\phi(x)\cdot y$$ for each $x\in E$. Clearly, $T$ is a positive compact operator (its rank is 1). But it is not an weak L-weakly compact operator.
	If not, as the singleton $\{u\}$ and
	$T(u) = \phi(u) · y = y$. For a norm bounded disjoint sequence $(y_n^\prime)\subset F^\prime$, clearly, $y^\prime_n(y)\nrightarrow0$, but $y_n^\prime\xrightarrow{un}0$. Therefore, $T$ is not weak L-weakly compact. It is absurd, so $F$ is order continuous, moreover $F$ is atomic. 
\end{proof}
Then we study that when is weak M-weakly compact operator compact.
\begin{theorem}\label{sc is compact}
	Let $E$ be a Dedekind $\sigma$-complete Banach lattice, then the following is equivalent.
	\begin{enumerate}
		\item $(E^\prime)^a$ is atomic and $E$ has order continuous norm.
		\item For each Banach space $Y$, every weak M-weakly compact operator $T:E\rightarrow Y$ is compact.
		\item For each Banach lattice $F$ without order continuous norm, every positive weak M-weakly compact operator $T:E\rightarrow F$ is compact.
	\end{enumerate}
\end{theorem}
\begin{proof}
	$(1)\Rightarrow(2)$. For a weakly M-weakly compact operator $T:E\rightarrow Y$. Since $E$ has order continuous norm, according to Theorem \ref{strong-M-weakly}, $T$ is M-weakly compact. It follows from \cite[Theorem~5.64]{AB:06} that $T^\prime:Y^\prime\rightarrow E^\prime$ is L-weakly compact.
	
	According to \cite[Proposition~3.6.2]{MN:91}, for any $\epsilon>0$, there exists some $x^\prime\in (E^\prime)^a_+$ such that $$T^\prime(B_{Y^\prime})\subset [-x^\prime,x^\prime]+\epsilon\cdot B_{E^\prime}.$$ Since $(E^\prime)^a$ is atomic, it follows from \cite[Theorem~6.1(5)]{W:99} that the order interval $[-x^\prime,x^\prime]$ is norm compact. Using \cite[Theorem~3.1]{AB:06}, we have $T(B_{Y^\prime})$ is relatively compact set in $E^\prime$, so $T^\prime$ is also compact. It follows Schauder theorem (\cite[Theorem~5.2]{AB:06}) that $T:E\rightarrow Y$ is compact.
	
	$(2)\Rightarrow(3)$. Obvious.
	
	$(3)\Rightarrow(1)$. First, we prove that $E$ has order continuous norm. Assume that $E$ is not order continuous, then $E$ and $F$ contain the lattice copy of $l_\infty$, moreover there exists a positive projection from $E$ to $l_\infty$. Let $R:E\rightarrow l_\infty$ be the positive projection, $S:l_\infty\rightarrow F$ be the canonical injection from $l_1$ into $F$ and $T=S\circ R:E\rightarrow l_\infty\rightarrow F$. Clearly, $T$ is a positive weak M-weakly compact operator, but not compact. Therefore, $E$ has order continuous norm.
	
	Then, we claim that $(E^\prime)^a$ is atomic. We assume that $(E^\prime)^a$ is not atomic and construct a
	positive weak M-weakly compact operator from $E$ into $F$ which is not compact.
	
	Since the norm of $F$ is not order continuous, According to  \cite[Theorem~4.14]{AB:06}, there
	exists a disjoint order bounded sequence $(y_n)$ in $F_+$ which does not converge to zero
	in norm. We may assume that $0 \leq y_n \leq y$ and $\Vert y_n\Vert=1$ for all $n$ and some $y\in F_+$. 
	
	Since $(E^\prime)^a$ is not atomic, there exists a $0\leq\psi\in (E^\prime)^a$ such that the order interval $[0,\psi]$ is not norm compact by \cite[Theorem~6.1]{W:99}. Choose a sequence $(\phi_n)$ in $[0,\psi]$ which
	has no norm convergent subsequence in $E^\prime$. Also, since $\psi\in (E^\prime)^a$, by \cite[Theorem~2.4.2]{MN:91}, the order interval
	$[0,\psi]$ is weakly compact. By the Eberlein-Smulian Theorem (\cite[Theorem~3.40]{AB:06}), we may assume, by extracting a subsequence if necessary, that $(\phi_n)$ converges
	weakly to some $\phi\in [0,\psi]$. So $(\phi_n)$ converges weakly* to $\phi$. 
	
	Now, define two operators $S,T : E \rightarrow F$ by
	$$S(x) = \phi (x) y+\sum_{n=1}^{\infty} (\phi_n-\phi)(x) y_n,\quad T(x)=\psi(x)y$$
	for each $x\in E$. It follows from the proof of Wickstead in \cite[Theorem~1]{W:96} that $0\leq S\leq T$ and that $S$ is not compact. 
	
	We claim that $S$ is weak M-weakly compact. For this, note that $0 \leq S^\prime \leq T^\prime$ and $T^\prime(h)=h(y)\psi$ for all $h\in F^\prime$. Then for every $h\in B_{F^\prime}$, we have
	$$|S^\prime(h)|\leq S^\prime(|h|)\leq T^\prime(|h|)\leq |h|y\psi\leq \Vert y\Vert\psi,$$
	Since $[-\psi,\psi]$ is L-weakly compact set, so $S^\prime$ is L-weakly compact operator, moreover $S^\prime$ is weak* L-weakly compact operator. According to Proposition \ref{weak L-strong continuous dual}, $S$ is weak M-weakly compact. Therefore, $(E^\prime)^a$ is atomic.
\end{proof}

Dually, we have:

\begin{corollary}\label{w*l is compact}
	Let $E$ be a Dedekind $\sigma$-complete Banach lattice, then the following statements are equivalent.
	\begin{enumerate}
		\item $(E^\prime)^a$ is atomic and $E$ has order continuous norm.
		\item For each Banach space $Y$, every weak* L-weakly compact adjoint operator $T^\prime:Y^\prime\rightarrow E^\prime$ for continuous operator $T:E\rightarrow Y$ is compact.
		\item For each Banach lattice $F$ without order continuous norm, every positive weak* L-weakly compact adjoint operator $T^\prime:F^\prime\rightarrow E^\prime$ for continuous operator $T: E\rightarrow F$ is compact.
	\end{enumerate}
\end{corollary}

The following result shows that when an operator is both compact and weak M-weakly compact.
\begin{theorem}\label{sc=compact}
	Let $E$ be a Dedekind $\sigma$-complete Banach lattice, then the following conditions are equivalent.
	\begin{enumerate}
		\item $E^\prime$ is atomic and both $E$ and $E^\prime$ are order continuous.
		\item For each Banach space $Y$, every continuous operator $T:E\rightarrow Y$ is weak M-weakly compact iff $T$ is compact.
		\item For each Banach lattice $F$ without order continuous norm, every positive operator $T:E\rightarrow F$ is weak M-weakly compact iff $T$ is compact.
	\end{enumerate}
\end{theorem}
\begin{proof}
	$(1)\Rightarrow(2)$. For a weakly M-weakly compact operator $T:E\rightarrow Y$. Since $E^\prime$ is atomic and $E$ is order continuous, it follows from Theorem \ref{sc is compact} that $T$ is compact. According to $E^\prime$ is order continuous and Schauder theorem (\cite[Theorem~5.2]{AB:06}), the converse is obtained immediately.
	
	$(2)\Rightarrow(3)$. Obvious.
	
	$(3)\Rightarrow(1)$ According to Theorem \ref{sc is compact}, we have $(E^\prime)^a$ is atomic and $E$ is order continuous. We claim that $E^\prime$ is order continuous.
	
	Assume that $E^\prime$ is not order continuous. There exists a positive projection $P:E\rightarrow l_1$. Fix a vector $0<y\in F_+.$ Define the operator $S:l_1\rightarrow F$ as follows:
	$$S\left(\lambda_{n}\right)=\left(\sum_{n=1}^{\infty} \lambda_{n}\right) y$$
	for each $(\lambda_n)\in l_1$. Obviously, the operator $S$ is well defined. Let
	$$T=S \circ P : E \rightarrow l_{1} \rightarrow F$$
	then $T$ is a positive compact operator since $S$ is a finite rank operator (rank is 1). Let $\left(e_{n}\right)$ be the standard basis of $l_{1}$. Clearly, $e_n\xrightarrow{un}0$, but $T(e_n)=y\nrightarrow0$. Hence, $T$ is not a weak M-weakly compact operator. This leads contradiction. Therefore $E^\prime$ is order continuous, moreover $E^\prime$ is atomic.
\end{proof}

The dual result is obtained immediately.

\begin{corollary}\label{w*l = compact}
	Let $E$ be a Dedekind $\sigma$-complete Banach lattice, then the following statements are equivalent.
	\begin{enumerate}
		\item $(E^\prime)^a$ is atomic and both $E$ and $E^\prime$ are order continuous.
		\item For each Banach space $Y$, every adjoint operator $T^\prime:Y^\prime\rightarrow E^\prime$ for continuous operator $T:E\rightarrow Y$ is weak* L-weakly compact iff $T$ is compact.
		\item For each Banach lattice $F$ without order continuous norm, every positive adjoint operator $T^\prime:F^\prime\rightarrow E^\prime$ for continuous operator $T:E\rightarrow F$ is weak* L-weakly compact iff $T$ is compact.
	\end{enumerate}
\end{corollary}

\noindent \textbf{Acknowledgement.} The research is supported by National Natural Science Foundation of China(NSFC:51875483).

\noindent \textbf{Data Availability Statement.} No data, models, or code were generated or used during the study.

\end{document}